\documentclass{notices}
\usepackage{amsfonts,amssymb,amsmath,amscd,xcolor} 
\usepackage[squaren]{SIunits}
\usepackage{graphicx}

\title{What is...Robin harmonic measure?}
\author{Max Engelstein
  \affil{
    Max Engelstein is an associate professor of mathematics at the University of Minnesota. His email address is mengelst@umn.edu.
    }
  \and
  Marcel Filoche
  \affil{
    Marcel Filoche is a CNRS Research Director at the Institut Langevin, ESPCI Paris. His email address is marcel.filoche@espci.psl.eu.
   }
   \and
   Svitlana Mayboroda
   \affil{
   Svitlana Mayboroda is a Professor of Mathematics at ETH Zurich. Her email address is svitlana.mayboroda@math.ethz.ch
   }}

\begin{document}

\maketitle

\begin{figure}
    \centering
    \includegraphics[width=0.45\linewidth]{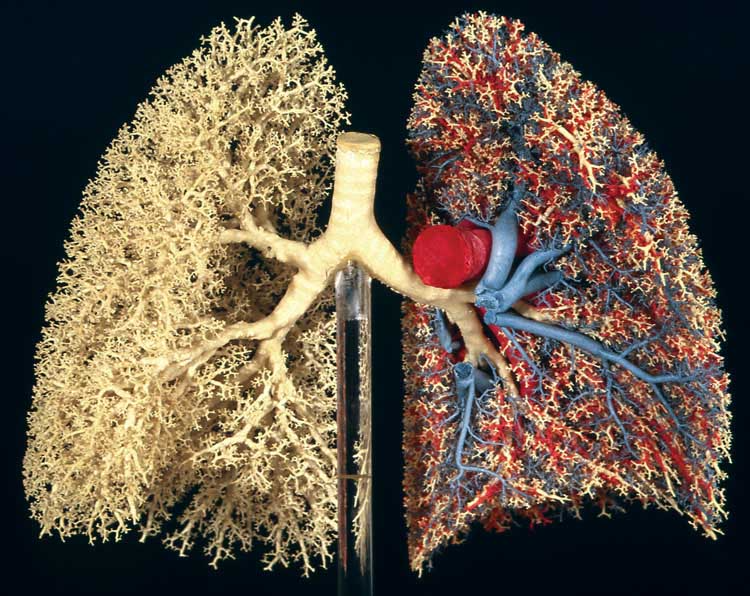}
    \hspace{2mm}
    \includegraphics[width=0.48\linewidth]{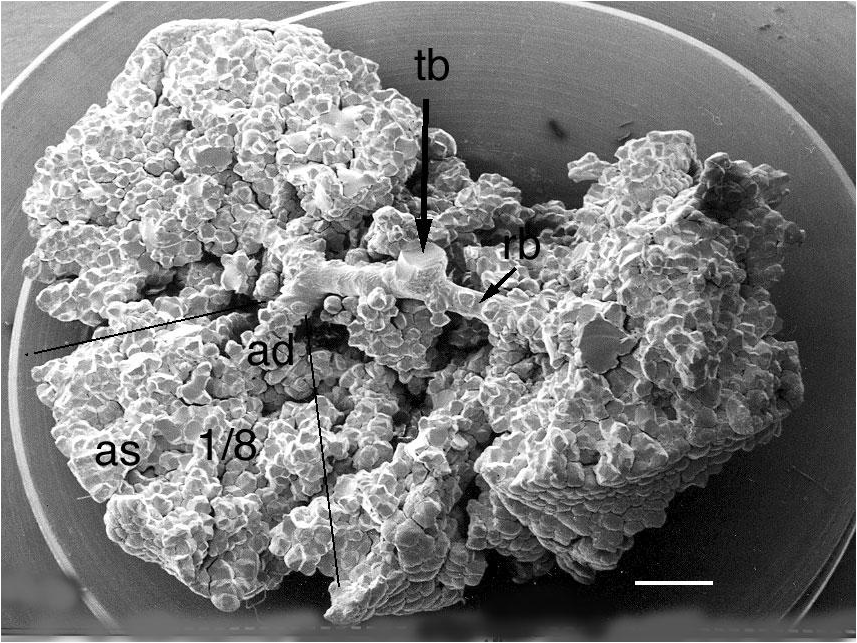}
    \caption{Left, Cast of a bronchial tree. Right, Human pulmonary acinus (scale bar: 1 mm). The acinus is the functional unit of gas exchange, where gas transport is dominated by diffusion rather than convection. (from~\cite{Haefeli-Bleuer_Weibel_1988})}
    \label{fig:lung_images}
\end{figure}

Diffusion-limited transport of gas or particles across membranes is ubiquitous in nature, especially in biological systems. A striking example is the human lung in which the inhaled air and the circulating blood are separated by the alveolo-capillary membrane through which oxygen and carbon dioxide counter-diffuse. Oxygen is advected by the air flow into the bronchial tree---a conducting tree of about 20 generations---then enters a region called \emph{acinus} where the ducts are connected to small chambers, the \emph{alveoli}~\cite{Haefeli-Bleuer_Weibel_1988}, see Fig.~\ref{fig:lung_images}. In the last five generations of the acinus (called \emph{subacinus}), oxygen transport is dominated by diffusion~\cite{Sapoval2002_PNAS}. The lung can therefore be seen as a collection of about 200,000 diffusion units operating in parallel.

Inside of one of these diffusion units, $\Omega$, let $u$ be the concentration of oxygen. The entrance of the unit (where it connects to the upper part of the bronchial tree) can be considered as the source of oxygen, whose concentration is normalized so that $u= 1$. Outside of this source but inside the diffusion unit it is reasonable to model the oxygen concentration as a solution to Laplace's equation. Along the wall of the diffusion unit the local oxygen flux can be assumed to be proportional to the concentration difference across the membrane (i.e. partial pressure), with the proportionality constant being determined by the membrane \emph{permeability}. These are known as Robin (or sometimes Fourier) boundary conditions. In summary we are interested in the solution, $u_a$ (to emphasize the dependence on~$a$) to the boundary value problem:
\begin{equation}\label{e:lungmodel}
\begin{cases} u_a = 1 \,,\quad &\textrm{in}\, B\subset \Omega,\\
 \Delta u_a = 0\,,\quad &\textrm{in}\, \Omega \backslash B,\\
 \displaystyle \frac{1}{a}\partial_\nu u_a + u_a = 0 \,,\quad &\textrm{in}\, \partial \Omega,
\end{cases}
\end{equation}
where $\nu$ is the outward pointing unit normal and $a$ is a constant proportional to the permeability of the lung membrane. If the permeability is infinite (perfect membrane), then $a=\infty$ and we recover a \emph{Dirichlet} boundary condition. On the contrary, if the permeability vanishes, then $a=0$ and we find a \emph{Neumann} boundary condition (blocking membrane).

The performance of the lung is related to its ability to transfer oxygen to the bloodstream. This is summarized by the total flux across the membrane
\begin{equation}
F(a):= -\int_{\partial \Omega}\partial_{\nu} u_a \, d\sigma \,,
\end{equation}
which depends on the absorption parameter~$a$ in a non-linear fashion. Furthermore, this non-linear dependence appears to be a function of the (pre-)fractal nature of $\Omega$.

It turns out that numerous pulmonary pathologies are related to either a modification of the domain or an alteration of its permeability. \emph{Pulmonary emphysema} consists of the progressive destruction of the alveolar walls, and hence a reduction of the alveolar surface available for gas exchange. In our framework, this corresponds to a reduction of the ``fractality'' of the boundary~$\partial \Omega$. \emph{Pulmonary edema} is an excessive accumulation of fluid in the air compartments of the lung. This fluid coats the alveolar surface, resulting in a reduction of membrane permeability, and ultimately impaired gas exchange. In our framework, this corresponds to a decrease of the absorption parameter~$a$. It is, for instance, one of the conditions associated with COVID-19.

Interestingly, it is well established in medicine that edema does not progress linearly with the amount of fluid present in the lung. Below a threshold of approximately \unit{300}{mL}, it is not a life-threatening condition: this is referred to as \emph{mild} edema. The lung is resilient. It is only above this threshold that shortness of breath appears even at rest, requiring either supplemental oxygen intake  or even mechanical ventilation or intubation. This is \emph{severe} edema. Understanding lung behavior and the existence of such a transition raises deep mathematical questions: (i)~How does the geometric complexity of the boundary, $\partial \Omega$, influence the value of the total flux,~$F(a)$? (ii)~How is gas transfer distributed over the boundary? (iii)~How is gas exchange affected by a change in the parameter~$a$? Answering these questions requires developing a precise understanding of the solution to~Eqs.~\eqref{e:lungmodel} at the boundary.

\subsection*{Classical (Dirichlet) harmonic measure: dimension drop}

We model the trajectory of an oxygen molecule inside a diffusion unit as a Brownian motion (a continuous random walk). When the molecule reaches the alveolar wall, it is absorbed into the bloodstream with probability $p$, and reflected back into the alveolus with probability $1-p$. The parameter $p$ is related to the absorptive constant $a$ by $a \approx \frac{p}{1-p}$. (The constant of proportionality is the diffusion mean free path)

In the idealized case of a perfectly absorbing boundary, i.e. when $p \equiv 1$ or equivalently $a=\infty$, the distribution describing where molecules first hit the boundary is called \emph{harmonic measure}. This is a classical object with remarkably subtle properties, and we begin by discussing this limiting case.

More precisely, let $\Omega \subset \mathbb R^n$ ($n\ge 2$) be a domain and let $x_0\in\Omega$. The harmonic measure $\omega_\Omega^{x_0}(E)$ of a set $E\subset \partial\Omega$ is the probability that a Brownian motion starting at $x_0$ first hits the boundary inside $E$ rather than $\partial\Omega\setminus E$.

Harmonic measure plays a central role in PDEs because of a theorem of Kakutani \cite{MR14646}, which represents solutions to the Dirichlet problem for the Laplacian in harmonic measure terms. Namely, for $f\in C(\partial\Omega)$, the function
$u_f(x):=\int_{\partial\Omega} f(Q)\, d\omega_\Omega^x(Q)$ where $x\in \Omega$
solves (for sufficiently regular $\Omega$)
\[
\begin{aligned}
\Delta u_f &= 0, \qquad && \textrm{in }\Omega,\\
u_f &= f, \qquad && \textrm{on }\partial\Omega.
\end{aligned}
\]

The relationship between the geometry of $\Omega$ and the behavior of $\omega_\Omega$ has been studied for more than a century, developing deep connections to singular integrals, geometric measure theory, and free boundary problems; see for example \cites{MR3540451, MR4169053, MR1726699}. Intuitively, if a Brownian particle must pass through a long narrow corridor, or repeatedly twist and turn, to reach a portion of the boundary $E\subset\partial\Omega$, then we expect $\omega_\Omega^{x_0}(E)$ to be small, even when $E$ itself is large. Some parts of the boundary may even be completely inaccessible, meaning they are hit with probability zero.

For fractal domains this phenomenon is often extreme: most of the boundary is typically inaccessible to Brownian motion. This is known as \emph{dimension drop}. Given a Radon measure $\mu$ on $\mathbb R^n$, define
\[
\dim_H(\mu)
:=
\inf\{\dim_H(E)\mid \mu(\mathbb R^n\setminus E)=0\},
\]
where $\dim_H$ denotes Hausdorff dimension. Roughly speaking, a set $E$ has Hausdorff dimension at least $\alpha$ if covering $E$ at scale $s\ll1$ requires on the order of $(1/s)^\alpha$ balls of radius $s$.

For smooth objects, geometric (Hausdorff) and harmonic measure dimensions coincide: a line has dimension $1$, a plane has dimension $2$, and harmonic measure is distributed along the entire boundary. Fractals behave very differently. For example, if $\Omega$ is bounded by the planar Koch snowflake, then
$\dim_H(\omega_\Omega^{x_0})=1,$
while $\dim_H(\partial\Omega)=\frac{\log 4}{\log 3}>1.$
Thus most of the boundary is effectively invisible to Brownian motion. In fact, the Koch snowflake is representative of a much broader phenomenon of dimension drop in fractal domains, see results in \cites{MR4140087, MR1404086}.

One might conjecture that $\dim_H(\omega_\Omega)\le n-1$ for every domain $\Omega\subset\mathbb R^n$. In two dimensions this is true, by celebrated results of Makarov \cite{MR794117} and Jones--Wolff \cite{MR962097}. In dimensions $n\ge3$, however, the picture changes dramatically. Wolff \cite{MR1315554} constructed simply connected domains $\Omega\subset\mathbb R^3$ for which $\dim_H(\omega_\Omega)>2,$ showing that the planar theory does not extend to higher dimensions. While Wolff's construction does not provide an explicit bound, numerical evidence suggests examples with $\dim_H(\omega_\Omega)>2.007$~\cite{GLFS2005}. These ideas even inspired the design of fractal noise-abatement walls installed along several roads in France \cite{FractalWall}.

Despite this, Bourgain proved a remarkable universal estimate \cite{MR874032}: for every dimension $n$ there exists $\epsilon_n>0$ such that $\dim_H(\omega_\Omega)\le n-\epsilon_n$ for all domains $\Omega\subset\mathbb R^n$. Determining the optimal value of $\epsilon_n$ remains one of the major open problems in the field. Recent work of Badger and Genschaw \cites{MR4887968, MR4824556} has refined Bourgain's argument and produced explicit quantitative bounds, although these are believed to be far from sharp.

 \subsection*{Robin Harmonic Measure and the Green function}

 In the previous section we referred to the Dirichlet harmonic measure, corresponding to the total absorption, as an {\it idealized} case. However, it turns out that it is not the optimal one in the sense of surface usage: it is the probability of reflection at the boundary that makes the Brownian motion explore all of the surface of the lung and thus ultimately makes it more efficient. Said otherwise, if the absorptive constant is infinite, there is very little incentive to develop more surface than a smooth, flat, $(n-1)$-dimensional one. It is the finite value of~$a$ that makes the development of a more complicated, rugged surface beneficial.

In \cite{david2024dimensionstructurerobinharmonic} the authors use PDE methods to rigorously construct the Robin Harmonic measure and establish its properties for a large class of domains $\Omega \subset \mathbb R^n$ (including many with fractal or pre-fractal boundaries). They show that there is a family of measures $\{\omega_{R, \Omega}^{x_0}\}_{x_0\in \Omega}$ such that for any $f\in C(\partial \Omega)$ the function $u_f^R(x_0): = \int_{\partial \Omega} f(Q) \, d\omega_{R,\Omega}^{x_0}(Q)$ satisfies:
 $$\begin{aligned}
     \Delta u_f^R =& 0\,\, \textrm{in } \Omega, \\
     \frac{1}{a}\partial_\nu u_f^R + u_f^R =& f\,\, \textrm{on } \partial \Omega.
 \end{aligned}$$ 

In stark contrast with the fully absorbing case, it has been shown in~\cite{david2024dimensionstructurerobinharmonic} that these measures {\bf never} demonstrate dimension drop. For all relevant $\Omega \subset \mathbb R^n$ and for all $a\in (0,\infty)$ 
$$\dim_H(\omega_{R, \Omega}^{x_0}) = \dim_H(\partial \Omega).$$
Heuristically, even just a tiny probability of reflection allows Brownian travelers to explore the full boundary. Returning to the lung model, \cite{david2024dimensionstructurerobinharmonic} implies that $u_a$ is H\"older continuous and there is a $\delta > 0$ (depending on $a$ and the geometry of $\Omega$) such that $$\inf_{\overline{\Omega}} u_a \geq \delta > 0.$$ That is, every part of the surface of the lung is transferring a quantitatively large amount of oxygen, thus partially justifying the lung's pre-fractal shape. All of it is used for the oxygen transfer, contrary to a seemingly more efficient totally absorbing case. 

From a physical perspective this result explains why the human lung has such a large surface, essentially exhibiting a fractal. It does not, however, explain the miraculous phase transition in the lung's performance that we have highlighted in the introduction. Here matters get technical, but let us at least give a taste of ideas.

The distribution of the harmonic measure along the boundary is intimately connected with the estimates on the Green function, the solution to 
$$ \begin{aligned}\Delta G^R(x, y) =\delta_y(x),\qquad & x\in \Omega,\\
     \frac{1}{a}\partial_\nu G^R(x, y) + G^R(x, y) = 0,\qquad &x\in \partial \Omega, \end{aligned}$$ 
for every $y\in \Omega$. In fact, on smooth domains, the Radon-Nykodim derivative of the harmonic measure (its distribution kernel) exists and is equal to the normal derivative of the Green function. This Robin Green function exhibits dramatically different behavior depending on the distance from $x$ to $\partial \Omega$. Morally speaking, the Green function is Neumann-like when $\mathrm{dist}(x, \partial \Omega)<1/a$ and Dirichlet-like when $\mathrm{dist}(x,\partial \Omega)>1/a$. We do not delve here into all the regimes connecting the distance between $x$ and $y$ and their distance to the boundary ---we refer the reader to \cite{david2025robingreenfunctionestimates} but it is this change from the Dirichlet to Neumann-like regime that defines the flux through the boundary and the phase transition that mystified us at the beginning.

\bibliography{biblio}

\end{document}